\documentclass{article}
\usepackage{camnum,psfig,harvard}
\newcommand{\gap}{\hspace*{-12.5pt}}

\begin{document}

\title{ Methods for the approximation of the matrix exponential
  in a Lie-algebraic setting
}

\author{Elena Celledoni\thanks{Research at MSRI is supported in part by NSF grant
DMS-9701755.},\\
MSRI,\\
1000 Centennial Drive,
Berkeley CA 94720,\\
celledon@msri.org,\\[0.2cm]
Arieh Iserles,\\
DAMTP,\\
Cambridge University,\\
Silver Street, CB3 9EW,\\
Cambridge, England\\
ai@damtp.cam.ac.uk.}

\maketitle

\begin{abstract}
  Discretization methods for ordinary differential equations based on
  the use of matrix exponentials have been known for decades. This set
  of ideas has come off age and acquired greater urgency recently,
  within the context of {\em geometric integration\/} and discretization
  methods on manifolds  based on the use of Lie-group actions.
  
  In the present paper we study the approximation of the matrix
  exponential in a particular context: given a Lie group $\CC{G}$ and
  its Lie algebra $\Gg{g}$, we seek approximants $F(tB)$ of $\exp(tB)$
  such that $F(tB)\in\CC{G}$ if $B\in \Gg{g}$.  Having fixed a basis
  $V_1,\dots, V_d$ of $\Gg{g}$, we write $F(tB)$ as a composition of
  exponentials of the type $\exp(\alpha_i(t)V_i)$, where $\alpha_i$ for
  $i=1,2,\ldots ,d$ are scalar functions.  In this manner it becomes
  possible to increase the order of the approximation without increasing
  the number of exponentials to evaluate and multiply together.  We
  study order conditions and implementation details and conclude the
  paper with some numerical experiments.
\end{abstract}

\section{Introduction}

Although numerical methods for the integration of ordinary
differential equations (ODEs) based on the use of the matrix
exponential have long history, the subject has acquired new relevance
recently with two developments. The first, which is irrelevant to the
theme of this paper, is the introduction of Krylov subspace techniques
and their application to large stiff systems of differential equations
\cite{hochbruck98eif}. The other development is motivated by the
philosophy of \textit{geometric integration\/} and its purpose is to
recover under discretization important qualitative and geometric
features of the underlying dynamical system. Examples of such methods
can be found {\em inter alia\/} in \cite{casas96ffa,crouch93nio}. An
important technique in geometric integration is  the use of Lie-group
actions, which lend themselves to the design of very effective
time-stepping methods for ODEs evolving on homogeneous manifolds. Such
methods have been recently studied in \cite{munthe-kaas97hor} and
\cite{engo98otc}. Methods based on the use of the classical Magnus and
Fer expansions for integrating ODEs on Lie-groups can be brought into
this formalism \cite{iserles97loi,iserles99ots,zanna97car}. All such methods
require a repeated evaluation of a matrix exponential, often of large
matrices. Inasmuch as typically one can expect the replacement of the
exact exponential by a suitable approximant (a rational function, say,
a Krylov subspace approximant or a Schur factorization), the context
of Lie-group methods imposes a crucial extra requirement. The
approximant in question, applied to an arbitrary element of the Lie
algebra $\GG{g}$, must produce an outcome in the Lie group $G$,
otherwise the whole purpose of the calculation, dicretizing within
$G$, will be null and void. This can be done is {\em some,\/} but by
no means, all Lie algebras of interest and we refer the reader to
\cite{celledoni98atm} for a more substantive discussion of this
issue. 

Let $G$ be a finite-dimensional Lie group. For all practical purposes,
we may assume that $G$ is a subgroup of the {\em general linear
  group\/} $\CC{GL}(n)$, the set of all nonsingular $n\times n$
matrices. We denote by $\GG{g}$ the Lie algebra corresponding to $G$,
observing that it is a subalgebra of $\GG{gl}(n)$, the Lie algebra of all
$n\times n$ matrices. Our concern in this paper is with differential
equations that evolve on a manifold ${\cal M}$ subject to the action
of $\CC{G}$. For simplicity we can assume that ${\cal M}$ coincides
with $\CC{G}$ and the action is of $\CC{G}$ on itself. The numerical
solution of such differential equations can be obtained considering
the {\em pull-back\/} on $\GG{g}$, by means of the exponential map,
of the vectorfield defining the equation.
We can compute the corresponding flow by a Lie-algebra discretization 
method and recover the approximation of the original problem via
exponentiation. 
Given an integration method of order $p$, we consider
order-$p$ approximants $F(tB)$ for $\exp (tB)$, where $B\in\GG{g}$ and
$t\geq0$. We require that $F(tB)\in \CC{G}$, whence it is easy to
prove that important qualitative features of the original equation and
the order of the discretization are retained. In \cite{celledoni98atm}
we have introduced low-rank splitting methods for the construction of
the approximant $F$, as the first attempt to provide a comprehensive
treatment of this issue.

Although the constraint $F(tB)\in \CC{G}$ represents remarkable
advantage in many applications, such as problems in which the
conservation of invariants is at issue in numerical modelling (volume
conservation in meteorology, invariance under rotations in the theory
of mechanical systems and in robotics), it should not be interpreted
as the sole purpose of our analysis. Our methods are relevant also for
the approximation of $\exp(tB)$ in the more general setting $B\in
\GG{gl}(n)$. Suppose in fact that $B\in \GG{gl}(n)$ and we want to
approximate $\exp(tB)$. It is always possible to write $B$ as a sum of
a matrix $B_s\in\GG{sl}(n)$ (the {\em special linear algebra\/} of
$n\times n$ matrices with zero trace) and a diagonal matrix $B_d$
whose nonzero entries are equal to $\delta=\CC{tr}\,(B)/n$. Then
$B_s=B-B_d$ and $[B_s ,B_d ]=0$ so that
$\exp(tB)=\exp(tB_s)\exp(t\delta)$. This fact is a particular case of
what is known in Lie theory as the Levi decomposition
\cite{humphreys72itl,varadarajan84lgl}.  Using this decomposition of
the matrix $B$, if necessary in tandem with some scaling and squaring
technique, the approximation of $\exp(tB)$ can be always reduced to
the approximation of $\exp(tB_s)$ with $B_s \in \GG{sl}(n)$. As long
as we can assure that our approximation of $\exp(tB_s)$ resides in
$\CC{SL}(n)$, the outcome is an approximant $F(tB)$ of $\exp(tB)$ that
shares with the exact exponential the feature that $\det F(tB)=\exp
(\CC{tr}\, B)$.

It is possible to prove that, given a splitting $B=\sum_{i=1}^{k}B_i$,
the function
\begin{displaymath}
  \CC{e}^{\frac12 tB_1}\CC{e}^{\frac12 tB_2}\cdots \CC{e}^{\frac12
  tB_{k-1}} \CC{e}^{tB_k}\CC{e}^{\frac12 tB_{k-1}}\cdots
  \CC{e}^{\frac12 tB_2} \CC{e}^{\frac12 tB_1},
\end{displaymath}
known as the generalized {\em Strang splitting,\/} approximates $\exp(tB)$ to
order 2. As long as $B_1,B_2,\ldots,B_k\in\GG{g}$, it follows at once
from the definition of a Lie group that the approximant resides in
$G$. Moreover, $2k-1$ is the least number of exponentials that
render such a splitting into a second-order approximant
\cite{celledoni98atm}. The Strang splitting is time reversible, hence
it follows readily from classical theory that the order can be raised
from 2 to 4 by composing three Strang splittings with different time
steps \cite{yoshida90coh}. In that case we need to evaluate $3k$
exponentials and multiply $6k$ matrices. 
In the case of low-rank splittings which have been considered by 
\citeasnoun{celledoni98atm} this results in the following count of
flops: $4n^3$ for order $2$, $12n^3$ for order
$4$. 

In this paper we present composition methods in which the number of
exponentials $k$ equals the dimension $d$ of the Lie algebra. Our
construction allows us to increase the order of the approximation
without increasing the number of exponentials to evaluate and multiply
together. Letting $\{V_1,\dots, V_d\}$ be a basis of $\GG{g}$, we
write $F(tB)$ as a composition of exponentials of the type
$\exp(\alpha_i(t)V_i)$, where each $\alpha_i(t)$ for $i=1,\dots ,d$ is
a scalar function.  In general $d=\O{ n^2}$, however, with an
appropriate choice of the basis elements, the computation of each
exponential $\exp(\alpha_i(t)V_i)$ requires $\O{n}$ flops, while the
formation of their product adds just $2n^3+\O{n^2}$ flops. The
challenging part of the computation is the construction of the
functions $\alpha_1(t) ,\dots ,\alpha_d(t)$, and the cost of their
calculation depends on the desired order of the approximation. Naive
complexity analysis might have indicated that the total cost is
growing exponentially in $d$ as the order increases. Yet, the cost
remains relatively modest for small orders and the method lends itself
very well to the exploitation of sparsity in the matrix $B$.  In the
sequel we show how this approach can be turned into an efficient
numerical method and we obtain algorithms of order up to $4$ with a
cost of $\O{n^3}$ for dense matrices.

Our approach can be interpreted as representing the solution using
{\em canonical coordinates of the second kind,\/} an approach that has
been pioneered by \citeasnoun{owren98imb} in the context of general
Lie-group methods. Having said this, the more restrictive framework of
exponential approximants possesses a very great deal of special
structure. This can be exploited so as to produce efficient and
competitive algorithms that approximate $\exp(tB)$, $B\in\GG{g}$, in
the Lie group $G$.

\section{The technique of coordinates of second kind for the
  approximation of the exponential matrix} 

Let $\CC{G}$ be a Lie group and $\GG{g}$ its corresponding
$d$-dimensional Lie algebra. We choose a basis $\{V_1,\dots, V_d\}$ of
$\GG{g}$, whence every element  $Y\in\CC{G}$ sufficiently close to the
identity can be represented in a unique fashion as
\begin{displaymath}
  Y=\exp(\gamma_1V_1) \exp(\gamma_2V_2)\cdots \exp(\gamma_dV_d),
\end{displaymath}
where $\exp : \GG{g}\rightarrow \CC{G}$ is the exponential map.  This
representation is known as representation in canonical coordinates of
the second kind \cite{varadarajan84lgl}. This representation is global
in the case of solvable Lie algebras.  We restrict ourselves to the
case $\GG{g}\subseteq\GG{gl}(n)$, $\CC{G}\subseteq\CC{GL}(n)$, when
$\exp$ is the usual matrix exponential.
 
Given $B\in\GG{g}$, we can represent it in a unique fashion as 
\begin{displaymath}
  B=\sum_{i=1}^{d}\beta_i V_i.
\end{displaymath}
It is possible then to write $\exp(tB)$ in the form
\begin{displaymath}
  U(t)=\exp(tB)=\exp(g_1(t)V_1) \exp(g_2(t)V_2)\cdots \exp(g_d(t)V_d).
\end{displaymath}
Letting $\MM{g}=[g_1,\dots,g_d]^{\CC{T}}$,
$\MM{\beta}=[\beta_1,\dots,\beta_d]^{\CC{T}},$  it can be proved that
the vector function $\MM{g}$ obeys a differential equation of the form
\begin{displaymath}
\frac{d\MM{g}}{dt}=\MM{f}(\MM{\beta},\MM{g}), \qquad \MM{g}(0)=0,
\end{displaymath}
where $\MM{f}$ is a suitable function of $\MM{\beta}$ and $\MM{g}$,
for sufficiently small $t$ \cite{wei63ogr}. Given a solvable Lie
algebra $\GG{g}$ \citeasnoun{wei63ogr} prove results on the global
representation of $U$. However, an explicit form of $\MM{f}$ is known
only for very simple examples of low-dimensional Lie algebras.  

In this paper we seek polynomials
${\alpha}_1\approx g_1,\dots,{\alpha}_d\approx g_d$ of a suitable
degree so that
\begin{displaymath}
  \exp(tB)\approx \exp({\alpha}_1(t)V_1) \exp({\alpha}_3(t)V_2) \cdots
  \exp({\alpha}_d(t)V_d). 
\end{displaymath}
Differentiation yields
\begin{equation}
  \label{eq:2.1}
  \sum_{i=1}^{d}\beta_i
  V_i=\sum_{i=1}^{d}g_i'(t)\prod_{j=1}^{i-1}\CC{e}^{g_jV_j}V_i
  \prod_{j=i-1}^{1}\CC{e}^{-g_jV_j}. 
\end{equation}
Evaluating this expression at the origin gives the first-order condition
\begin{equation}
  \label{eq:2.2}
  g'_i(0)=\beta_i,\qquad i=1,2,\ldots,d.
\end{equation}
Further differentiations of \R{eq:2.1} lead to  higher-order
conditions.  Let us define the functions 
\begin{equation}
  \label{eq:2.3}
  P_i(\MM{g})=\exp(\CC{ad}_{g_1V_1})\circ
  \cdots\circ\exp(\CC{ad}_{g_{i-1}V_{i-1}})(V_i),\qquad i=1,2,\ldots,d,
\end{equation}
where the {\em adjoint operator\/} $\CC{ad}_{x}:\GG{g}\rightarrow
\GG{g}$ is defined as $\CC{ad}_{x}(y)=[x,y]$ for any $x,y \in\GG{g}$,
$[x,y]=xy-yx$ being the matrix commutator. Note that
$P_i(\MM{g}(0))=V_i$, $i=1,2,\ldots,d$. Moreover, the right-hand side
of  \R{eq:2.1} can be written in the simplified form
The function
\begin{displaymath}
  {\cal T}(\MM{g})=\sum_{i=1}^{d}g_i'(t)P_i(g).
\end{displaymath}
Since the derivatives of the left hand side of \R{eq:2.1} vanish, the
conditions for order $p\geq 1$, can be obtained by solving the equations 
\begin{equation}
  \label{eq:2.4}
  \left.\frac{\CC{d}^r}{\CC{d}t^r}{\cal T}(\MM{g})\right|_{\,t=0}=0,
  \qquad r=1,2,\ldots , p-1, \quad p\geq 1,
\end{equation}
where
\begin{equation}
  \label{eq:2.5}
  \frac{\CC{d}^r}{\CC{d}t^r}{\cal
  T}(\MM{g})=\sum_{i=1}^{d}\sum_{k=1}^{r}{r\choose k}
  \frac{\CC{d}^{r-k+1}g_i}
  {\CC{d}t^{r-k+1}}\frac{\CC{d}^{k}P_i}{\CC{d}t^{k}}. 
\end{equation}
In particular,
\begin{Eqnarray*}
  \frac{\CC{d}}{\CC{d}t}{\cal T}(\MM{g})&=&\sum_{i=1}^{d} \left(g_i''P_i
    +g_i'\frac{d}{dt}P_i\right),\\
  \frac{\CC{d}^2}{\CC{d}t^2}{\cal T}(\MM{g})&=&\sum_{i=1}^{d} \left(
    g_i'''P_i+2g_i''\frac{\CC{d}}{\CC{d}t}P_i+g_i'\frac{\CC{d}^2}
    {\CC{d}t^2}P_i\right) ,\\
  \frac{\CC{d}^3}{\CC{d}t^3}{\cal T}(\MM{g})&=&\sum_{i=1}^{d} \left(
    g_i^{IV}P_i+3g_i'''\frac{\CC{d}}{\CC{d}t}P_i+3g_i''
    \frac{\CC{d}^2}{\CC{d}t^2}P_i+g_i'\frac{\CC{d}^3}{\CC{d}t^3}P_i\right).
\end{Eqnarray*}

Solving  \R{eq:2.4} for $r=1$ results in the values of $g_i''(0)$ for
$i=1,2,\ldots ,d$ that allow us to construct an order-$2$ 
approximant. Substituting such values in \R{eq:2.4} for $r=2$ yields
$g_i'''(0)$ for $i=1,2,\ldots ,d$ and consequently an approximant of 
order $3$. Similar procedure can be used to construct recursively
approximants of arbitrarily high order.

The main part of the computation is the evaluation of the $k$-th
derivative of $P_i(\MM{g})$ at $t=0$. Expanding the exponentials in
\R{eq:2.3} we obtain 
\begin{displaymath}
  P_i(\MM{g})=\prod_{k=1}^{i-1} (I+\CC{ad}_{g_kV_k}+\Frac12
  \CC{ad}^2_{g_kV_k}+\Frac16 \CC{ad}^3_{g_kV_k}+\ldots )(V_i)
\end{displaymath}
and, after further algebra,
\begin{Eqnarray*}
  P_i(\MM{g})&=&\left\{I+\sum_{k=1}^{i-1}\CC{ad}_{g_kV_k}+
    \sum_{k=2}^{i-1}
    \sum_{l=1}^{k-1}\CC{ad}_{g_lV_l}\CC{ad}_{g_kV_k}\right.\\
  &&\mbox{}+\Frac12\sum_{k=1}^{i-1}\CC{ad}^2_{g_kV_k}+
  \sum_{k=3}^{i-1}\sum_{l=2}^{k-1}\sum_{j=1}^{l-1}\CC{ad}_{g_jV_j}
  \CC{ad}_{g_lV_l}\CC{ad}_{g_kV_k}\\
  &&\left.\mbox{}+\Frac12\sum_{k=2}^{i-1}\sum_{l=1}^{k-1}
  \left(\CC{ad}_{g_lV_l}\CC{ad}^2_{g_kV_k}+\CC{ad}^2_{g_lV_l}
    \CC{ad}_{g_kV_k}\right)+\Frac16\sum_{k=1}^{i-1}\CC{ad}^3_{g_k
    V_k}+\dots \right\}\left(V_i\right).
\end{Eqnarray*}
Similarly to \cite{owren97rkm}, we write $P_i(\MM{g})$ in the
form
\begin{displaymath}
  P_i(\MM{g})=I+\sum_{r=1}^{\infty}\sum_{j_1=1}^{i-1}
  \sum_{j_2=j_1}^{i-1}\cdots \sum_{j_r=j_{r-1}}^{i-1}
  \frac{1}{\MM{j}!}g_{j_1}\dots
  g_{j_r}\CC{ad}_{V_{j_1}}\circ\cdots\circ\CC{ad}_{V_{j_r}}(V_i). 
\end{displaymath}
Here $\MM{j}=(j_1,\dots j_r)$ is a multi-index of integer elements with
$1\leq j_r\le i-1$ and $\MM{j}!:=q_1!q_2!\dots q_{i-1}!$ where $q_k$ is the
number of occurrences of $k$ in $(j_1,j_2\ldots j_r)$.

A general expression for the $k$-th derivative of $P_i$ is given
as follows: since $g_i(0)=0$, we may let $f_i(t)=g_i(t)/t$,
$i=1,2,\ldots,d$. We can then rewrite $P_i$ in the form
\begin{displaymath}
  P_i(\MM{g})=
  I+\sum_{r=1}^{\infty}t^r\sum_{j_1=1}^{i-1}\sum_{j_2=j_1}^{i-1}\cdots
  \sum_{j_r=j_{r-1}}^{i-1}\frac{1}{\MM{j}!}f_{j_1}\cdots
  f_{j_r}\CC{ad}_{V_{j_1}}\circ\cdots\circ\CC{ad}_{V_{j_r}}(V_i).
\end{displaymath}
By following the construction in\cite{owren97rkm} we obtain
\begin{equation}
  \label{eq:2.6}
  \begin{meqn}
    \displaystyle \left. \frac{\CC{d}^kP_i}{\CC{d}t^k}\right|_{t=0}&=&
    \displaystyle \sum_{r=1}^{k} \sum_{\delta_1+\dots+\delta_r=k}
    \frac{k!}{\prod_{\nu =1}^{\mu}(\delta_{\nu}-1)!}\sum_{1\le
      j_1\le\dots\le j_{\mu}\le
      i-1}\frac{1}{\MM{j}!}\\[18pt]
    &&\left.\displaystyle \mbox{}\times f_{j_1}^{(\delta_1-1)}\cdots
      f_{j_{\mu}}^{(\delta_{\mu}-1)}\right|_{t=0}
    \CC{ad}_{V_{j_1}}\circ\cdots\circ \CC{ad}_{V_{j_{\mu}}}(V_i).  
  \end{meqn}
\end{equation}

Substituting \R{eq:2.6} in \R{eq:2.4} and \R{eq:2.5} we obtain the conditions
for arbitrary order $p$. In particular we obtain the following
formulae for the derivatives of $P_i(\MM{g})$ at $t=0$,
\begin{Eqnarray*}
  \left.\frac{\CC{d}P_i}{\CC{d}t}\right|_{t=0}&=&\sum_{k=1}^{i-1}
  \CC{ad}_{V_k}\left(V_i\right)g_k'(0),\\ 
  \left. \frac{\CC{d}^2P_i}{\CC{d}t^2}\right|_{t=0}&=&\sum_{k=1}^{i-1}
  \left(\sum_{l=1}^{k-1}2\CC{ad}_{V_l}\CC{ad}_{V_k}(V_i)g_k'(0)
    g_l'(0)+\CC{ad}^2_{V_k}(V_i)[g_k'(0)]^2+\CC{ad}_{V_k}(V_i)g_k''(0)\right),\\
  \left. \frac{\CC{d}^3P_i}{\CC{d}t^3}\right|_{t=0}&=&
  6\sum_{k=3}^{i-1}\sum_{l=2}^{k-1}\sum_{j=1}^{l-1}\CC{ad}_{V_j}
  \CC{ad}_{V_l}\CC{ad}_{V_k}(V_i)g_j'(0)g_l'(0)g_k'(0)\\ 
  &&\mbox{}+3\sum_{k=2}^{i-1}\sum_{l=1}^{k-1}\left(\CC{ad}_{V_l}
    \CC{ad}^2_{V_k}(V_i)g_l'(0)[g_k'(0)]^2+\CC{ad}^2_{V_l}\CC{ad}_{V_k}
    (V_i) [g_l'(0)]^2g_k'(0)\right)\\
  &&\mbox{}+\sum_{k=1}^{i-1}\CC{ad}_{V_k}^3(V_i)[g_k'(0)]^3\\
  &&\mbox{}+3\sum_{k=2}^{i-1}\sum_{l=1}^{k-1}\CC{ad}_{V_l}
  \CC{ad}_{V_k} (V_i)\left(g_k''(0)g_l'(0)+g_k'(0)g_l''(0)\right)\\
  &&\mbox{}+3\sum_{k=1}^{i-1}\CC{ad}^2_{V_k}(V_i)g_k''(0)g_k'(0)\\
  &&\mbox{}+\sum_{k=1}^{i-1}\CC{ad}_{V_k}(V_i)g_k'''(0).
\end{Eqnarray*}
Substitution readily produces order conditions. Specifically, 
\begin{equation}
  \label{eq:2.7}
  \sum_{i=1}^{d}
  g_i''(0)V_i=-\sum_{i=1}^{d}g_i'(0)\sum_{k=1}^{i-1}\CC{ad}_{V_k}(V_i)g_k'(0),
\end{equation}
are conditions for order $2$, while
\begin{Eqnarray}  
  \sum_{i=1}^{d}  g_i'''(0)V_i&=& -  \sum_{i=1}^{d}\left\{
    2g_i''(0)\sum_{k=1}^{i-1}\CC{ad}_{V_k}(V_i)g_k'(0)\right. \nonumber\\
  &&\mbox{}+g_i'(0)\sum_{k=1}^{i-1}\left[2\sum_{l=1}^{k-1}
    \CC{ad}_{V_l}\CC{ad}_{V_k}(V_i)g_k'(0)g_l'(0)\right.\label{eq:2.8}\\
  &&\left.\left.\mbox{}+\CC{ad}^2_{V_k}(V_i)(g_k'(0))^2+\CC{ad}_{V_k}
    (V_i)g_k''(0)\right]\right\},\nonumber 
\end{Eqnarray}
are the order-3 conditions. Finally, conditions for order 4 are
\begin{Eqnarray}
  \sum_{i=1}^{d} g_i^{\CC{IV}}(0)V_i&=&- \sum_{i=1}^{d}\left\{
    3g_i'''(0)\sum_{k=1}^{i-1}\CC{ad}_{V_k}(V_i)g_k'(0)+\right. \nonumber\\ 
  &&\mbox{}+3g_i''(0)\left[\sum_{k=1}^{i-1}\left(\sum_{l=1}^{k-1}
      2\CC{ad}_{V_l}\CC{ad}_{V_k}(V_i)g_k'(0)g_l'(0)\right.\right.\nonumber\\
  &&\left.\left.\mbox{}+ \CC{ad}^2_{V_k}(V_i)(g_k'(0))^2+\CC{ad}_{V_k}
      (V_i)g_k''(0)\right)\right]\nonumber\\
  &&\mbox{}+g_i'(0)\left[6\sum_{k=3}^{i-1}\sum_{l=2}^{k-1}\sum_{j=1}^{l-1}
    \CC{ad}_{V_j}\CC{ad}_{V_l}\CC{ad}_{V_k}(V_i)g_j'(0)g_l'(0)g_k'(0)\right. 
  \nonumber\\ 
  &&\mbox{}+3\sum_{k=2}^{i-1}\sum_{l=1}^{k-1}\left(\CC{ad}_{V_l}
    \CC{ad}^2_{V_k}(V_i)g_l'(0)(g_k'(0))^2+\CC{ad}^2_{V_l}\CC{ad}_{V_k} 
    (V_i)(g_l'(0))^2g_k'(0)\right)\label{eq:2.9}\\
  &&\mbox{}+\sum_{k=1}^{i-1}\CC{ad}_{V_k}^3(V_i)(g_k'(0))^3\nonumber\\
  &&\mbox{}+3\sum_{k=2}^{i-1}\sum_{l=1}^{k-1}\CC{ad}_{V_l}\CC{ad}_{V_k} 
  (V_i)\left(g_k''(0)g_l'(0)+g_k'(0)g_l''(0)\right)\nonumber\\
  &&\mbox{}+3\sum_{k=1}^{i-1}\CC{ad}^2_{V_k}(V_i)g_k''(0)g_k'(0)
  \nonumber\\ 
  &&\mbox{}+\left.\left. \sum_{k=1}^{i-1}\CC{ad}_{V_k}(V_i)g_k'''(0)
    \right]\right\}. \nonumber
\end{Eqnarray}

\begin{figure}[tb]
  \begin{center}
    \leavevmode
    \psfig{file=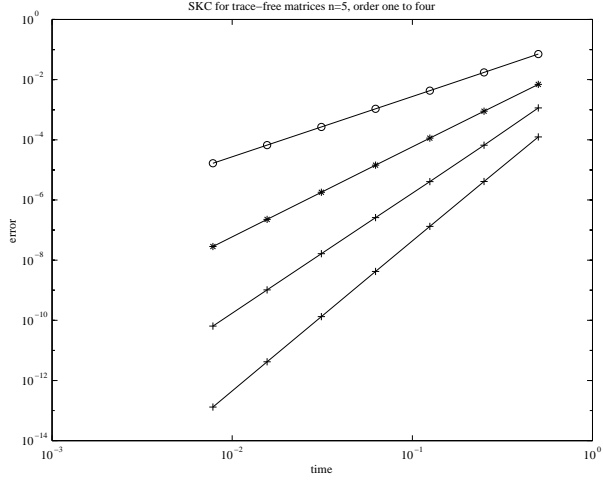,width=8cm}
    \caption{Error in the approximation of the exponential with WN
      technique.}
    \label{fig:1}
  \end{center}
\end{figure}

In Figure \ref{fig:1} we have plotted along the $y$ axis the
$2$-norm of the error of the approximation of $\exp(tB)$ with the
second-kind coordinates (SKC) methods of order ranging from $1$ to
$4$. The values of the error are plotted against time, (along the
$x$-axis), to logarithmic scale for matrices of $\GG{sl}(5)$. The
methods have been implemented using the standard basis defined in
section 3.

The computation of $g''(0)$, $g'''(0)$ and $g^{\CC{IV}}(0)$ is
obtained directly implementing the formulas \R{eq:2.7},
\R{eq:2.8},\R{eq:2.9} respectively. This implementation does not
depend on the choice of the particular basis of $\GG{sl}(5)$, but the
number of commutators that must be computed with this approach is
$\O{d^p}$ for $p=2,3,4$. Even if we assume that the $V_i$s are very
sparse matrices and that the cost of computing each commutator is
$\O{1}$ operations, the total cost exceeds $\O{n^{2p}}$ flops for
$p=2,3,4$ where $n$ is the dimension of the matrix. Such expense is
not acceptable for a competitive method of approximation of
$\exp(tB)$. Fortunately, it can be decreased a very great deal by an
appropriate choice of the basis $\{V_1,V_2,\ldots,V_d\}$. This is the
theme of the next section.

\setcounter{equation}{0}
\section{Choosing a basis}

The choice of the right basis and sparse representation of commutators
are critical to the implementation of the SKC methods. Recalling the
order conditions \R{eq:2.7}--\R{eq:2.9}, our aim is to choose a basis
so that terms of the form $\CC{ad}_{V_{i_1}}\CC{ad}_{V_{i_2}} \cdots
\CC{ad}_{V_{i_s}} V_j$ can be represented in the most economical
manner. We recall that, given the basis $\{V_1,V_2,\ldots,V_d\}$ of a
$d$-dimensional Lie algebra $\GG{g}$, the {\em structure constants\/}
are the numbers $c_{k,l}^i$, $k,l,i=1,2,\ldots,d$, such that
\begin{displaymath}
  [V_k,V_l]=\sum_{i=1}^d c_{k,l}^i V_i
\end{displaymath}
\cite{humphreys72itl}. Let 
\begin{displaymath}
  B=\sum_{k=1}^d \beta_k V_k.
\end{displaymath}
Then an order-1 condition is always
\begin{equation}
  \label{eq:3.1}
  g_k'(0)=\beta_k,\qquad k=1,2,\ldots,d.
\end{equation}
To obtain the order-2 condition we substitute \R{eq:3.1} in \R{eq:2.7}
and express commutators in terms of structure constants,
\begin{Eqnarray*}
  \sum_{k=1}^d g_k''(0)V_k&=&-\sum_{l=1}^d \beta_l \sum_{j=1}^{l-1}
  [V_j,V_l] \beta_j=-\sum_{l=1}^d \beta_l \sum_{j=1}^{l-1}
  \sum_{k=1}^d c_{j,l}^k V_k\\
  &=&-\sum_{k=1}^d \left( \sum_{l=1}^d \sum_{j=1}^{l-1} \beta_l
  c_{j,l}^k \beta_j \right) V_k.
\end{Eqnarray*}
Since $c_{j,l}^k=-c_{l,j}^k$, we thus deduce that
\begin{equation}
  \label{eq:3.2}
  g_k''(0)=\sum_{l=1}^d \sum_{j=1}^{l-1} \beta_l c_{l,j}^k
  \beta_j,\qquad k=1,2,\ldots,d. 
\end{equation}
Likewise, substituting in \R{eq:2.8},
\begin{Eqnarray*}
  \sum_{k=1}^d g_k'''(0)V_k&=&-\sum_{i=1}^d \left\{
    2g_i''(0)\sum_{l=1}^{i-1} [V_l,V_i]\beta_l +\beta_i \sum_{l=1}^{i-1}
    \left(2\sum_{j=1}^{l-1} [V_j,[V_l,V_i]] \beta_l\beta_j
    \right.\right.\\
  &&\left.\left.\mbox{}+[V_l,[V_l,V_i]]\beta_l^2
      +[V_l,V_i]g_l''(0)\right) \right\}.
\end{Eqnarray*}
Note that
\begin{displaymath}
  [V_j,[V_l,V_i]]=\sum_{s=1}^d c_{l,i}^s [V_j,V_s] =\sum_{k=1}^d
  \sum_{s=1}^d c_{l,i}^s c_{j,s}^k V_k.
\end{displaymath}
Therefore
\begin{Eqnarray*}
  \sum_{k=1}^d g_k'''(0)V_k&=&-2\sum_{i=1}^d g_i''(0) \sum_{l=1}^{i-1}
  \sum_{k=1}^d c_{l,i}^k \beta_l V_k-2\sum_{i=1}^d \beta_i
  \sum_{l=1}^{i-1} \sum_{j=1}^{l-1} \sum_{k=1}^d \sum_{s=1}^d
  c_{l,i}^s c_{j,s}^k \beta_l \beta_j V_k\\
  &&\mbox{}-\sum_{i=1}^d \beta_i \sum_{l=1}^{i-1} \sum_{k=1}^d
  \sum_{s=1}^d c_{l,i}^s c_{l,s}^k \beta_l^2 V_k-\sum_{i=1}^d
  \sum_{l=1}^{i-1} \beta_i \sum_{k=1}^d c_{l,i}^k g_l''(0)V_k
\end{Eqnarray*}
and we deduce that
\begin{Eqnarray*}
  g_k'''(0)&=&\sum_{i=1}^d \sum_{l=1}^{i-1} c_{i,l}^k
  [2g_i''(0)\beta_l +\beta_i g_l''(0)]+ 2\sum_{i=1}^d \sum_{l=1}^{i-1}
  \sum_{j=1}^{l-1} \sum_{s=1}^d c_{i,l}^s c_{j,s}^k \beta_i\beta_l \beta_j\\
  &&\mbox{}-\sum_{i=1}^d\sum_{i=1}^d \sum_{l=1}^{i-1} \sum_{s=1}^d
  c_{i,l}^s c_{l,s}^k \beta_i \beta_l^2,\qquad k=1,2,\ldots,d.
\end{Eqnarray*}
Bearing in mind that for order $p$ we require
\begin{displaymath}
  \alpha_k(t)=\sum_{r=1}^p \frac{1}{r!} g_k^{(r)}(0)t^r,\qquad
  k=1,2,\ldots,d, 
\end{displaymath}
we observe that the sheer volume of calculations required for the
evaluation of the functions $\alpha_1,\alpha_2,\ldots,\alpha_d$ is
prohibitive for, say, order 3, unless most of the structure constants
vanish.  Fortunately, bases of finite-dimensional Lie algebras which
are `sparse' (in the sense that a very high proportion of structure
constants vanish) are known. They are associated with {\em root space
decompositions\/} of Lie algebras \cite{humphreys72itl} and, in the
case of semisimple algebras, are known as {\em Chevalley bases\/}
\cite{carter95llg}. Wishing to avoid too much Lie-algebraic
terminology in a numerical analysis paper, we reserve our exposition
to just three examples which are the most important in a range applications.

\vspace{8pt}
\noindent {\bf The orthogonal group} 
Let $\GG{g}=\GG{so}(n)$, the Lie algebra of $n\times n$ skew-symmetric
matrices. It corresponds to two important Lie groups: the orthogonal
group $\CC{O}(n)$ of $n\times n$ orthogonal matrices and its subgroup,
the special orthogonal group $\CC{SO}(n)$ of matrices with unit
determinant. Its dimension is $d=\frac12n(n-1)$. We let
\begin{displaymath}
  F_{i,j}=\MM{e}_i{\MM{e}_j}^{\CC{T}}-\MM{e}_j{\MM{e}_i}^{\CC{T}},
  \qquad i=1,2,\ldots,n, \quad j=i+1,i+2,\ldots ,n,
\end{displaymath} 
where $\MM{e}_i$ is the $i$-th canonical vector of $\BB{R}^n$. In
other words, $F_{i,j}$ is a matrix whose $(i,j)$-th element is $1$, the
$(j,i)$-th element equals $-1$ and zero otherwise.
We can trivially expand each $B\in \GG{so}(n)$ as
$B=\sum_{i=1}^{n}\sum_{j=i+1}^nb_{i,j}F_{i,j}$. $U(t):=\exp(tF_{i,j})$
is simply an {\em Euler rotation\/} in the $(i,j)$ plane: it is
identity matrix, except that 
\begin{displaymath}
  \left[
    \begin{array}{cc}
      U_{i,i} & U_{i,j}\\
      U_{j,i} & U_{j,j}
    \end{array}\right]=\left[
    \begin{array}{rr}
      \cos(tF_{i,j}) & \sin(tF_{i,j})\\
      -\sin(tF_{i,j}) & \cos(tF_{i,j}) 
    \end{array}\right].
\end{displaymath}
Noting that
\begin{displaymath}
  [F_{i,j} , F_{l,k} ] = 
  \begin{case}
    -F_{j,k}, &  i=l,\; j\neq k,\\
    -F_{i,l}, & i\neq l,\; j=k, \\
    F_{i,k}, & i\neq k,\; j=l, \\
    F_{j,l}, & i=k,\; j\neq l, \\
    O, & \CC{otherwise,}
  \end{case}
\end{displaymath}
the order conditions are simplified as follows,
\begin{Eqnarray*}
      p\ge 1:\qquad g'_{i,j}(0)&=&b_{i,j},\\
      p\ge 2:\qquad g''_{i,j}(0)&=&
      \sum_{s=j+1}^nb_{j,s}b_{i,s}-\sum_{r=i+1}^{j-1}b_{r,j}b_{i,r}\\
      &&\mbox{}+\sum_{r=1}^{i-1}b_{r,j}b_{r,i},\qquad 1\leq i<j\leq n,
\end{Eqnarray*}
and similarly for higher-order terms. Thus, the cost of computing the
coefficients for the second-order method is just $\frac12(n-2)(n-1)n
\approx\frac12 n^3$ flops. In comparison, a naive computation of
\R{eq:3.2}, without exploiting sparsity of structure constants,
requires $\frac18 (n^2-n-2)(n-1)^2n^2\approx \frac18 n^6$ flops. 

A more classical composition method for $B\in\GG{so}(n)$ is the
\textit{Strang splitting\/} which we can write in the form
\begin{displaymath}
  \CC{e}^{tb_{1,2}F_{1,2}/2}\cdots
  \CC{e}^{tb_{n-2,n}F_{n-2,n}/2} \CC{e}^{tb_{n-1,n}F{n-1,n}}
  \CC{e}^{tb_{n-2,n}F_{n-2,n}/2}\cdots 
  \CC{e}^{tb_{1,2}F_{1,2}/2}
\end{displaymath}
\cite{celledoni98atm}. It gives a second-order approximant to
$\exp(tB)$ whose calculation requires $\approx 4n^3$ flops, in
comparison with $\approx 3n^3$ for the second-order CSK method.

We note that, in the specific case of $\GG{so}(n)$, diagonal
Pad\'e approximants to the exponential provide an alternative to our
method, since they map the algebra to $\CC{O}(n)$. Having said this,
for dense matrices $B$ the cost of evaluating the second-order
approximant $(I-\frac12 tB)^{-1}(I+\frac12 tB)$ with, say, LU
factorization is $\O{n^3}$, comparative with our method.

\vspace{8pt}
\noindent {\bf The special linear group}
Let $\GG{g}=\GG{sl}(n)$, the set of $n\times n$ matrices with zero
trace, whence $d=n^2-1$. We split the  algebra in the first instance
into  diagonal and off-diagonal parts: in the terminology of Lie
algebras, the subspace spanned by the diagonal elements is a {\em
  Cartan subalgebra\/} \cite{carter95llg} or {\em maximal toral
  algebra\/} \cite{humphreys72itl} of $\GG{sl}(n)$. Specifically, our
basis is
\begin{displaymath}
  \{E_{i,j}\,:\, i,j=1,2,\ldots,n\; i\neq j\}\cup \{D_i\,:\,
  i=1,2,\ldots,n-1\}. 
\end{displaymath}
where
\begin{Eqnarray*}
  E_{i,j}&=&\MM{e}_i\MM{e_j}^{\CC{T}}, \qquad i,j=1,2,\ldots , n,\quad
  i\neq j,\\
  D_i&=&\MM{e}_i\MM{e_i}^{\CC{T}}-\MM{e}_{i+1}\MM{e}_{i+1}^{\CC{T}},
  \qquad i=1,2,\ldots,n-1.  
\end{Eqnarray*}
The exponentials of $E_{i,j}$ and $D_i$ are trivial,
\begin{displaymath}
  \CC{e}^{tE_{i,j}}=I+tE_{i,j},\qquad \CC{e}^{tD_i}=\CC{e}^t D_i.
\end{displaymath}
We order the elements by taking first $E_{i,j}$, $i\neq j$, in
lexicographic order, followed by $D_1,D_2,\ldots,D_{n-1}$. The
commutator table is
\begin{Eqnarray*}
  [E_{i,j},E_{r,s}]&=&
  \begin{case}
    E_{i,s}, & i\neq s,\; j=r,\\
    -E_{r,j}, & i=s,\;j\neq r,\\
    \sum_{l=s}^{r-1} D_l, & i=s<j=r,\\
    -\sum_{l=r}^{s-1} D_l, & i=s>j=r,\\
    O, & \mbox{otherwise,}
  \end{case}\qquad 
  \begin{array}{l}
    i,j,r,s=1,2,\ldots,n,\\
    i\neq j,\quad r\neq  s,
  \end{array} \\{} 
  [E_{i,j},D_r]&=&
  \begin{case}
    -E_{r,j}, & i=r,\; j\neq r+1,\\
    -E_{i,r+1}, & i\neq r,\; j=r+1,\\
    -2E_{r,r+1}, & i=r,\; j=r+1,\\
    E_{r+1,j}, & i=r+1,\; j\neq r,\\
    E_{i,r}, & i\neq r+1,\; j=r,\\
    2E_{r+1,r}, & i=r+1,\; j=r,\\
    O, & \mbox{otherwise,}
  \end{case}\qquad 
  \begin{array}{l}
    i,j=1,\ldots,n,\quad i\neq j,\\
    r=1,\ldots,n-1,
  \end{array}\\{}
  [D_i,D_j]&=&O,\qquad i,j=1,2,\ldots,n-1.
\end{Eqnarray*}
In general, for a $d$-dimensional Lie algebra there are $(d-1)d^2$
structure constants. In the case of $\GG{sl}(n)$ this means that up to
$\approx n^6$ structure constants may be nonzero. Yet, using the above
basis results in just $2(n-1)n^2+4(n-2)(n-1)+\frac23 (n^2-1)n\approx
\frac73 n^3$ nonzero structure constants and substantial saving in the
implementation of the SKC technique. 

Letting
\begin{Eqnarray*}
  [E_{i,j},E_{r,s}]&=&\sum_{(k,l)} c_{(i,j),(k,l)}^{(k,l)} E_{k,l} +
  \sum_k c_{(i,j),(r,s)}^{(k)} D_k,\\{}
  [E_{i,j},D_r]&=&\sum_{(k,l)}c_{(i,j),r}E_{k,l}
\end{Eqnarray*}
(note that $[D_r,E_{i,j}]=-[E_{i,j},D_r]$ and $[D_i,D_j]=O$) we thus
have 
\begin{Eqnarray*}
  c_{(i,j),(r,s)}^{(k,l)}&=&
  \begin{case}
    +1, & k=i,\; l=s,\; r=j,\; s\neq i,\\
    -1, & k=r,\; l=j,\; r\neq j,\; s=i,\\
    0, & \mbox{otherwise},
  \end{case}\\
  c_{(i,j),(r,s)}^{k}&=&
  \begin{case}
    +1, & i=s<j=r,\; k=s,s+1,\ldots,r-1,\\
    -1, & i=s>j=r,\; k=r,r+1,\ldots,s-1,\\
    0, & \mbox{otherwise},
  \end{case}\\
  c_{(i,j),r}^{(k,l)}&=&
  \begin{case}
    +1, & k=i=r+1,\;l=j,\;j\neq r \mbox{\ or\ } k=i,\; l=j=r,\;
    i\neq r+1,\\
    -1, & k=i=r,\;l=j,\;j\neq r+1 \mbox{\ or\ } k=i,\; l=j=r+1,\;
    i\neq r,\\
    +2, & k=i=r+1,\; l=j=r,\\
    -2, & k=i=r,\; l=j=r+1,\\
    0, & \mbox{otherwise},
  \end{case}\\
  c_{(i,j),r}^{k}&=&c_{r,s}^{(k,l)}=c_{r,s}^k=0.
\end{Eqnarray*}
Letting
\begin{displaymath}
  B=\sum_{k\neq l} \beta_{k,l} E_{k,l}+\sum_{k} \gamma_k D_k,
\end{displaymath}
and ordering the pairs $(k,l)$, $k\neq l$, in lexicographic order, we
thus have 
\begin{Eqnarray*}
  g_{k,l}'(0)&=&\beta_{k,l},\\
  g_k'(0)&=&\gamma_k,\\
  g_{k,l}''(0)&=&\sum_{(i,j)\succ(r,s)} \beta_{i,j}
  c_{(i,j),(r,s)}^{(k,l)} \beta_{r,s} +\sum_{(i,j),r} \beta_{i,j}
  c_{(i,j),r}^{k} \gamma_r\\
  &=&\sum_{i=1}^{k-1} \beta_{k,i}\beta_{i,l}-\sum_{i=k+1}^n
  \beta_{k,i}\beta_{i,l} +\beta_{k,l}(\gamma_{k-1}+\gamma_l -\gamma_k
  -\gamma_{l-1}),\\
  g_k''(0)&=&\sum_{(i,j)\succ(r,s)} \beta_{i,j} c_{(i,j),(r,s)}^{k}
  \beta_{r,s} =-\sum_{i=1}^k \sum_{j=k+1}^n \beta_{i,j}\beta_{j,i},
\end{Eqnarray*}
where $\gamma_0=\gamma_n=0$.

\vspace{8pt}
\noindent{\bf The Lorenz group} This is the 6-dimensional group
$\CC{SO}(3,1)$ of $4\times4$ matrices $A$ such that $AJA^{\CC{T}}=J$,
where $J=\CC{diag}(1,1,1,-1)$ \cite{carter95llg}. It has important
applications in special relativity theory.  he corresponding {\em
  Lorenz algebra\/} $\GG{so}(3,1)$ consists of all matrices $B$ such
that $BJ+JB^{\CC{T}}=O$. It is easy to verify that each element of
$\GG{so}(3,1)$ can be written in the form
\begin{displaymath}
  B=\left[
    \begin{array}{rrrr}
      0 & b_1 & b_2 & b_3\\
      -b_1 & 0 & b_4 & b_5\\
      -b_2 & -b_4 & 0 & b_6\\
      b_3 & b_5 & b_6 & 0
    \end{array}\right],\qquad b_1,b_2,\ldots,b_6\in\BB{R}. 
\end{displaymath}
Choosing the basis
\begin{Eqnarray*}
  &&\left\{\left[
      \begin{array}{rrrr}
        0 & 1 & 0 & 0\\
        -1 & 0 & 0 & 0\\
        0 & 0 & 0 & 0\\
        0 & 0 & 0 & 0
      \end{array}\right],\;\, \left[
      \begin{array}{rrrr}
        0 & 0 & 1 & 0\\
        0 & 0 & 0 & 0\\
        -1 & 0 & 0 & 0\\
        0 & 0 & 0 & 0
      \end{array}\right],\;\, \left[
      \begin{array}{rrrr}
        0 & 0 & 0 & 0\\
        0 & 0 & 1 & 0\\
        0 & -1 & 0 & 0\\
        0 & 0 & 0 & 0
      \end{array}\right],\;\, 
    \left[
      \begin{array}{rrrr}
        0 & 0 & 0 & 1\\
        0 & 0 & 0 & 0\\
        0 & 0 & 0 & 0\\
        1 & 0 & 0 & 0
      \end{array}\right],\;\, 
    \left[
      \begin{array}{rrrr}
        0 & 0 & 0 & 0\\
        0 & 0 & 0 & 1\\
        0 & 0 & 0 & 0\\
        0 & 1 & 0 & 0
      \end{array}\right],\;\, \left[
      \begin{array}{rrrr}
        0 & 0 & 0 & 0\\
        0 & 0 & 0 & 0\\
        0 & 0 & 0 & 1\\
        0 & 0 & 1 & 0
      \end{array}\right]\right\},
\end{Eqnarray*}
we obtain the commutator table
\begin{displaymath}
  \begin{array}{llllll}
    \gap&[V_1,V_2]=-V_3, &[V_1,V_3]=V_2, &[V_1,V_4]=-V_5, &[V_1,V_5]=V_4,
    &[V_1,V_6]=O,\\{} 
    \gap&[V_2,V_1]=V_3, &[V_2,V_3]=-V_1, &[V_2,V_4]=-V_6, &[V_2,V_5]=O,
    &[V_2,V_6]=V_4,\\{} 
    \gap&[V_3,V_1]=-V_2, &[V_3,V_2]=V_1, &[V_3,V_4]=O, &[V_3,V_5]=-V_6,
    &[V_3,V_6]=V_5,\\{}
    \gap&[V_4,V_1]=V_5, &[V_4,V_2]=V_6, &[V_4,V_3]=O, &[V_4,V_5]=V_1,
    &[V_4,V_6]=V_2,\\{} 
    \gap&[V_5,V_1]=-V_4, &[V_5,V_2]=O, &[V_5,V_3]=V_6,
    &[V_5,V_4]=-V_1, &[V_5,V_6]=V_3,\\{} 
    \gap&[V_6,V_1]=O, &[V_6,V_2]=-V_4, &[V_6,V_3]=-V_5,
    &[V_6,V_4]=-V_2, &[V_6,V_5]=-V_3.
  \end{array}
\end{displaymath}
Thus, out of 180 structure constants, just 24 are nonzero -- and they
all equal $\pm1$. After brief claculation, we drive for example the
polynomials $\alpha_k$ that yield an order-2 CSK approximant,
\begin{Eqnarray*}
  \alpha_1(t)&=&\beta_1 t+\Frac12(\beta_2\beta_3-\beta_4\beta_5)t^2,\\
  \alpha_2(t)&=&\beta_2 t-\Frac12(\beta_1\beta_3+\beta_4\beta_6)t^2,\\
  \alpha_3(t)&=&\beta_3 t+\Frac12(\beta_1\beta_2-\beta_5\beta_6)t^2,\\
  \alpha_4(t)&=&\beta_4 t-\Frac12(\beta_1\beta_5+\beta_2\beta_6)t^2,\\
  \alpha_5(t)&=&\beta_5 t+\Frac12(\beta_1\beta_4-\beta_3\beta_6)t^2,\\
  \alpha_6(t)&=&\beta_6 t+\Frac12(\beta_3\beta_5-\beta_2\beta_4)t^2,
\end{Eqnarray*}
where $B=\sum_{k=1}^6 \beta_k V_k$.

\setcounter{equation}{0}
\section{Time symmetry}

An approximant $F(tB)\approx\exp(tB)$ is said to be {\em time
  symmetric\/} if $F(tB)F(-tB)=I$, $t\geq0$. Time symmetric
approximants are important for a number of reasons, not least being
that they lend themselves to the {\em Yo\v{s}ida technique,\/} which
allows their order to be increased \cite{yoshida90coh}. The techniques
of the last section are not time symmetric. Here we describe their
modification, which results in a time-symmetric approximant.

Me mention in passing that it is possible to envisage two distinct
techniques to obtain high-order algorithms based on canonical
coordinates of the second kind. The first, implicit in the work of the
previous section, consists of evaluating the numbers $g_k^{(l)}(0)$
for $l=1,2,\ldots,p$, where $p$ is the order of the method. The
alternative, the subject matter of the present section, consists in
combining a second-order or a fourth-order approximant across a number
of steps to obtain a higher-order method.

Given the splitting
\begin{displaymath}
  B=\sum_{l=1}^s C_l,
\end{displaymath}
it is well known that the \emph{Strang splitting\/}
The approximation
\begin{equation}
  \label{eq:4.1}
  F(tB)=\CC{e}^{tC_1/2}\cdots \CC{e}^{tC_{s-1}/2}e^{tC_s}
  \CC{e}^{tC_{s-1}/2}\cdots \CC{e}^{tC_1/2} 
\end{equation}
is of order $2$ and time symmetric. Note that, as a consequence of
time symmetry, for sufficiently small $t\geq0$ we can represent
$F(tB)=\CC{e}^{{\cal F}(t)}$ where the matrix function ${\cal F}(t)$
is odd. It is precisely this feature that allows the application of
the Yo\v{s}ida technique.

The clear reason for \R{eq:4.1} being time symmetric is that it is
palindromic in the alphabet $\{C_1,C_2,\ldots,C_s\}$. This provides a clue how to
modify techniques based on canonical coordinates of the second kind so
as to render them time symmetric. Given a basis
$\{V_1,V_2,\ldots,V_d\}$ of the Lie algebra $\GG{g}$, we approximate
$\CC{e}^{tB}$ by the product
\begin{equation}
  \label{eq:4.2}
  \exp[\alpha_1(t)V_1]\cdots \exp[\alpha_{d-1}(t)V_{d-1}]
  \exp[\alpha_d(t)V_d] \exp[\alpha_{d-1}(t)V_{d-1}]\cdots
  \exp[\alpha_{1}(t)V_1], 
\end{equation}
where $\alpha_1,\alpha_2,\ldots,\alpha_d$ are {\em odd\/} polynomials.

Taking $\alpha_l=\frac12\beta_l t$, $l=1,2,\ldots,d-1$ and
$\alpha_d=\beta_d t$ yeilds the second-order Strang splitting. In the
sequel we seek higher-order methods of this kind.

Using the Baker--Campbell--Hausdorff (BCH) formula 
it is possible to express the product of exponentials at the right
hand side of \R{eq:4.2} as a single exponential \cite[
p. 141]{varadarajan84lgl}. Due to the symmetric arrangement of the
exponentials in \R{eq:4.2}, the BCH formula is an expansion in odd
powers of $t$. If this expansion converges, which is always the case
for sufficently small $t$, it makes sense to write the equation 
\begin{equation}
  \label{eq:4.3}
  tB=2\sum_{i=1}^{d-1}\alpha_i(t)V_i+\alpha_d(t)V_d+
  \sum_{k=1}^{\infty}Q^{2k}(\MM{\alpha}).  
\end{equation}
Here we denote by $Q^{2k}(\MM{\alpha})$ the terms of order
$\O{t^{2k+1}}$ in the  BCH formula applied to \R{eq:4.3}.
Moreover, we let $\alpha_{i}^{2k}(t)$ be the polynomial
obtained by truncating  the expansion of $\alpha_i(t)$ after the first
$k$ terms, and we denote the remainder by $r_i^{2k}(t)$. In
other words,
\begin{displaymath}
  \alpha_i(t)=\alpha_i^{2k}(t)+r_i^{2k}(t),\qquad
  r_i^{2k}(t)=\O{t^{2k+1}},\qquad i=1,2,\ldots,d.
\end{displaymath}
From \R{eq:4.3} we deduce
\begin{displaymath}
  2\sum_{i=1}^{d-1}\alpha_i^{2(k)}(t)V_i+\alpha_d^{2(k)}(t)V_d=tB-
  \sum_{r=1}^{k-1}Q^{2r}(\MM{\alpha})+\O{t^{2k+1}}.  
\end{displaymath}
Noting that 
\begin{displaymath}
  Q^{2r}(\MM{\alpha})=Q^{2r}(\MM{\alpha}^{2(k-1)}+\MM{r}^{2(k-1)})
  =Q^{2r}(\MM{\alpha}^{2(k-1)})+\O{t^{2k+r}}, 
\end{displaymath}
we obtain
\begin{equation}
  \label{eq:4.4}
  2\sum_{i=1}^{d-1}\alpha_i^{2(k)}(t)V_i+\alpha_d^{2(k)}(t)V_d
  =tB-\sum_{r=1}^{k-1}Q^{2r}(\MM{\alpha}^{2(k-1)})+\O{t^{2k+1}}. 
\end{equation}
Dropping the $\O{t^{2k+1}}$ terms in \R{eq:4.4}, it is possible to compute
$\MM{\alpha}^{2k}$ from $\MM{\alpha}^{2(k-1)}$. This gives a procedure  to derive a
sequence of successively increaing-order approximants of $\exp(tB)$.
It is easy to see that the approximants
\begin{displaymath}
  F^{2k}(tB)=\exp(\alpha_1^{2k}(t)V_1)\cdots
  \exp(\alpha_d^{2k}(t)V_d)\cdots \exp(\alpha_{1}^{2k}(t)V_1), 
\end{displaymath}
of $\exp(tB)$ are such that $F^{2k}(tB)F^{2k}(-tB)=I$, hence time
symmetry, the reason being the symmetric arrangements of the
exponentials in $F^{2k}(tB)$ and the odd-power expansion of the
functions $\alpha_i^{2k}$.

The BCH and symmetric BCH formulae for $k$-terms have an exceedingly 
complicated expansion, which can be obtained recursively. In what
follows we will make use just of the term $Q^{2}(\alpha)$,
demonstrating how it is possible to compute  it explicitely for
particular choices of the basis.

In the remainder of this section we consider the implementation of
time-symmetric CSK methods. We split $B$ as before and commence by
considering the Strang splitting \R{eq:4.1} except that, to simplify
notation, we arrange the terms in reverse ordering,
\begin{equation}
  \label{eq:4.5}
  \CC{e}^{tC_s/2}\cdots \CC{e}^{tC_2/2}\CC{e}^{tC_1}\CC{e}^{tC_2/2}\cdots
  \CC{e}^{tC_s/2}. 
\end{equation}

\begin{lemma}
  The term $Q^2$ of the BCH formula applied to \R{eq:4.5} is
  \begin{equation}
    \label{eq:4.6}
    Q^2= \frac{t^3}{12}\sum_{l=2}^s [C_1+\cdots+C_{l-1} +\Frac12
    C_l,[C_1+\cdots+C_{l-1},C_l]]. 
  \end{equation}
\end{lemma}
\begin{proof}
  See the appendix.
\end{proof}

Let us next consider the case $\GG{g}=\GG{so}(n)$, choosing the same
sparse basis as in Section~2. Therefore, according to \R{eq:4.6}, we
have
\begin{equation}
  \label{eq:4.7}
  \begin{meqn}
    Q^2&=&\displaystyle \frac{t^3}{12}\sum_{i=1}^{n-1}\sum_{j=i+1}^n
    [b_{1,2}F_{1,2}+ \ldots +
    b_{i,j-1}F_{i,j-1},[b_{1,2}F_{1,2}+\ldots+b_{i,j-1} 
    F_{i,j-1},\\
    &&\mbox{}\displaystyle b_{i,j}F_{i,j}]]+
    \frac{1}{24}\sum_{i=1}^{n-1}\sum_{j=i+1}^n  b_{i,j} 
    [F_{i,j},[b_{1,2}F_{1,2}+\dots+b_{i,j-1}F_{i,j-1},b_{i,j}F_{i,j}]].
  \end{meqn}
\end{equation}
We compute separately each part of this sum. Exploiting the commutator
table of our basis we have  
\begin{Eqnarray*}
  &&[b_{1,2}F_{1,2}+\ldots+b_{i,j-1}F_{i,j-1},b_{i,j}F_{i,j}]\\
  &=& b_{i,j}\left(-\sum_{s=i+1}^{j-1}b_{i,s}F_{s,j}-\sum_{r=1}^{i-1}
  b_{r,j}F_{r,i}+\sum_{t=1}^{i-1}b_{t,i}F_{t,j}\right), 
\end{Eqnarray*}
and noting that $b_{i,s}=-b_{s,i}$, we deduce that
\begin{equation}
  \label{eq:4.8}
  [b_{1,2}F_{1,2}+\dots+b_{i,j-1}F_{i,j-1},b_{i,j}F_{i,j}]=b_{i,j}\left(
    \sum_{\stackrel{\scriptstyle t=1}{t\neq
    i}}^{j-1}b_{t,i}F_{t,j}-\sum_{r=1}^{i-1}b_{r,j}F_{r,i}\right). 
\end{equation}
Commuting the right-hand side with $F_{i,j}$ gives
\begin{equation}
  \label{eq:4.9}
  \left[F_{i,j},\sum_{\stackrel{\scriptstyle t=1}{t\neq
  i}}^{j-1}b_{t,i}F_{t,j}-\sum_{r=1}^{i-1}
  b_{r,j}F_{r,i}\right]= \sum_{\stackrel{\scriptstyle t=1}{t\neq
  i}}^{j-1} b_{t,i}F_{t,i}+\sum_{r=1}^{i-1}b_{r,j}F_{r,j}. 
\end{equation}
Let $\MM{b}_1,\MM{b}_2,\ldots,\MM{b}_n$ be the columns of $B$ and
denote 
\begin{displaymath}
  \MM{b}_l^{s}=\sum_{k=1}^{s-1}b_{k,l}\MM{e}_k,\qquad k=1,2,\ldots,n.
\end{displaymath}
Then  \R{eq:4.8} yields
\begin{displaymath}
  [b_{1,2}F_{1,2}+\dots+b_{i,j-1}F_{i,j-1},b_{i,j}F_{i,j}]=b_{i,j}(
    \MM{b}_{i}^{j}\MM{e}_j^{\CC{T}}-\MM{e}_j{\MM{b}_{i}^{j}}^{\CC{T}})
    -b_{i,j}(\MM{b}_{j}^{i}\MM{e}_i^{\CC{T}}
    -\MM{e}_i{\MM{b}_{j}^{i}}^{\CC{T}}), 
\end{displaymath}
while \R{eq:4.9} gives
\begin{Eqnarray*}
  &&[F_{i,j}, [b_{1,2}F_{1,2}+\cdots
  +b_{i,j-1}F_{i,j-1},b_{i,j}F_{i,j}]]\\
  &=& b_{i,j}(\MM{b}_{i}^{j}\MM{e}_i^{\CC{T}}-\MM{e}_i{\MM{b}_{i}^{j}}
  ^{\CC{T}})-b_{i,j}(\MM{b}_{j}^{i}\MM{e}_j^{\CC{T}}-\MM{e}_j
  {\MM{b}_{j}^{i}}^{\CC{T}}).
\end{Eqnarray*}
Multiplying the latter by $b_{i,j}$ and summing in $i$ and $j$ we can
evaluate  \R{eq:4.7} in $n^3$ operations. Note that we count
separately multiplications and additions, for example, we assume that
the cost of Euclidean inner product of two vectors of length $n$ is 
$2n$ operations.

We now assemble together our results to calculate \R{eq:4.7}. We
proceed by splitting the sum
$b_{1,2}F_{1,2}+\cdots+b_{i,j-1}F_{i,j-1}$ in three parts, whereby 
\begin{Eqnarray*}
  &&[b_{1,2}F_{1,2}+\cdots+b_{i,j-1}F_{i,j-1},\MM{b}_{i}^{j}
  \MM{e}_j^{\CC{T}}-\MM{e}_j{\MM{b}_{i}^{j}}^{\CC{T}}-(\MM{b}_{j}^{i}
  \MM{e}_i^{\CC{T}} -\MM{e}_i{\MM{b}_{j}^{i}}^{\CC{T}})]\\
  &=&\sum_{l=1}^{i}\sum_{k=1}^{l-1}[(\MM{b}_{l}^{l}\MM{e}_l^{\CC{T}}
  -\MM{e}_l{\MM{b}_{l}^{l}}^{\CC{T}}),\MM{b}_{i}^{j}\MM{e}_j^{\CC{T}}
  -\MM{e}_j{\MM{b}_{i}^{j}}^{\CC{T}}-(\MM{b}_{j}^{i}\MM{e}_i^{\CC{T}} 
  -\MM{e}_i{\MM{b}_{j}^{i}}^{\CC{T}})]\\
  &&\mbox{}+\sum_{l=i+1}^{j-1}\sum_{k=1}^{i}[(\MM{b}_{l}^{i+1}\MM{e}_l^{\CC{T}}
  -\MM{e}_l{\MM{b}_{l}^{i+1}}^{\CC{T}}),
  \MM{b}_{i}^{j}\MM{e}_j^{\CC{T}} -\MM{e}_j{\MM{b}_{i}^{j}}^{\CC{T}}
  -(\MM{b}_{j}^{i}\MM{e}_i^{\CC{T}}-\MM{e}_i{\MM{b}_{j}^{i}}^{\CC{T}})]\\
  &&\mbox{}+\sum_{l=j}^{m}\sum_{k=1}^{i-1}[(\MM{b}_{l}^{i}\MM{e}_l^{\CC{T}}
  -\MM{e}_l{\MM{b}_{l}^{i}}^{\CC{T}}), \MM{b}_{i}^{j}\MM{e}_j^{\CC{T}}
  -\MM{e}_j{\MM{b}_{i}^{j}}^{\CC{T}}-(\MM{b}_{j}^{i}\MM{e}_i^{\CC{T}}
  -\MM{e}_i{\MM{b}_{j}^{i}}^{\CC{T}})].
\end{Eqnarray*}
Finally,
\begin{Eqnarray*}
  &&[b_{1,2}F_{1,2}+\cdots+b_{i,j-1}F_{i,j-1},\MM{b}_{i}^{j}
  \MM{e}_j^{\CC{T}}-\MM{e}_j{\MM{b}_{i}^{j}}^{\CC{T}}-(\MM{b}_{j}^{i}
  \MM{e}_i^{\CC{T}} -\MM{e}_i{\MM{b}_{j}^{i}}^{\CC{T}})]\\
  &=&-{\MM{b}_i^i}^{\CC{T}}\MM{b}_i^jF_{i,j}+\MM{b}_i^i
  {\MM{b}_j^i}^{\CC{T}}-\MM{b}_j^i{\MM{b}_i^i}^{\CC{T}}+
  {\MM{b}_j^{i}}^{\CC{T}}\MM{b}_j^iF_{j,i}-(\MM{b}_j^i
  {\MM{b}_i^j}^{\CC{T}}-\MM{b}_i^j{\MM{b}_j^i}^{\CC{T}})\\   
  &&\mbox{}\sum_{l=1}^{i-1}b_{l,i}(\MM{b}_l^l\MM{e}_j^{\CC{T}}-
  \MM{e}_j{\MM{b}_l^l}^{\CC{T}})-{\MM{b}_l^l}^{\CC{T}}\MM{b}_i^jF_{l,j} 
  -b_{l,j}(\MM{b}_l^l\MM{e}_i^{\CC{T}}-\MM{e}_i{\MM{b}_l^l}^{\CC{T}})
  +{\MM{b}_l^l}^{\CC{T}}\MM{b}_j^iF_{l,i}\\ 
  &&\mbox{}\sum_{l=i+1}^{j-1}b_{l,i}(\MM{b}_l^{i+1}\MM{e}_j^{\CC{T}}
  -\MM{e}_j{\MM{b}_l^{i+1}}^{\CC{T}})-{\MM{b}_l^{i+1}}^{\CC{T}}
  \MM{b}_i^{j}F_{l,j}+b_{i,l}(\MM{b}_{j}^{i}\MM{e}_l^{\CC{T}}-
  \MM{e}_l{\MM{b}_{j}^{i}}^{\CC{T}})+{\MM{b}_l^{i+1}}^{\CC{T}}\MM{b}_j^i F_{l,i}\\ 
  &&\mbox{}+\sum_{l=j+1}^{n}{\MM{b}_l^{i}}^{\CC{T}}\MM{b}_j^iF_{l,i}
  -{\MM{b}_l^i}^{\CC{T}}\MM{b}_i^jF_{l,j}. 
\end{Eqnarray*}
We analyse the computational costs of the previous formula, summing
over $i$ and $j$ and showing that \R{eq:4.7} can be computed in
$\O{n^3}$ operations. Note that, since
$\sum_{l=i+1}^{j-1}b_{i,l}\MM{e}_l=\MM{b}_i^i-\MM{b}_i^j$, we have
\begin{displaymath}
  \MM{b}_i^i{\MM{b}_j^i}^{\CC{T}}-\MM{b}_j^i{\MM{b}_i^i}^{\CC{T}}+
  \sum_{l=i+1}^{j-1}b_{i,l}(\MM{b}_{j}^{i}\MM{e}_l^{\CC{T}}-
  \MM{e}_l{\MM{b}_{j}^{i}}^{\CC{T}})-(\MM{b}_j^i{\MM{b}_i^j}^{\CC{T}}
  -\MM{b}_i^j{\MM{b}_j^i}^{\CC{T}})=-2(\MM{b}_j^i{\MM{b}_i^j}^{\CC{T}}
  -\MM{b}_i^j{\MM{b}_j^i}^{\CC{T}}).
\end{displaymath}
It is more convenient to write the previous expression in the form
\begin{Eqnarray*}
  &&-2\sum_{i=1}^{n-1}\sum_{j=i+1}^n b_{i,j}(\MM{b}_j^i
  {\MM{b}_i^j}^{\CC{T}}-\MM{b}_i^j{\MM{b}_j^i}^{\CC{T}})\\
  &=&-\sum_{i=1}^n\left( \sum_{j=i+1}^{n}2b_{i,j}\MM{b}_j^i\right)
  {\MM{b}_i^i}^{\CC{T}}-{\MM{b}_i^i}\left(\sum_{j=i+1}^{n}2b_{i,j}
    \MM{b}_j^{i}\right)^{\!\!\CC{T}}\\  
  &&\mbox{}-\sum_{i=1}^{n-1}\sum_{j=i+1}^n2b_{i,j}
  \left(\MM{b}_j^i(\MM{b}_i^j-\MM{b}_i^i)^{\CC{T}}
  -(\MM{b}_i^j-\MM{b}_i^i){\MM{b}_j^i}^{\CC{T}}\right).  
\end{Eqnarray*}
The first part of this sum is computed in about $\frac23 n^3$ operations and the
second part, exploiting the equality
\begin{displaymath}
  \sum_{i=1}^{n-1}\sum_{j=i+1}^n2b_{i,j}\MM{b}_j^i(\MM{b}_i^j
  -\MM{b}_i^i)^{\CC{T}}=2\sum_{i=1}^{n-1}\sum_{k=n-1}^{i+2}b_{i,k}
  \left(\sum_{l=n}^{k+1}b_{i,l}\MM{b}_l^i\right)\MM{e}_k^{\CC{T}}, 
\end{displaymath}
can also be computed in $\frac23 n^3$ operations.

Adding terms of the type $\alpha F_{l,j}$ and $\beta F_{l,i}$ leads to
\begin{displaymath}
  -\sum_{l=1}^i{\MM{b}_l^l}^{\CC{T}}\MM{b}_i^jF_{l,j}-
  \sum_{l=i+1}^{j-1}{\MM{b}_l^{i+1}}^{\CC{T}}\MM{b}_i^{j}F_{l,j}
  -\sum_{l=j+1}^{n}{\MM{b}_l^i}^{\CC{T}}\MM{b}_i^jF_{l,j}
  =-L_i^i\MM{b}_j^i{\MM{e}_j}^{\CC{T}}-{\MM{e}_j}(-L_i^i\MM{b}_j^j)^{\CC{T}},   
\end{displaymath}
and
\begin{displaymath}
  \sum_{l=1}^{i-1}{\MM{b}_l^l}^{\CC{T}}\MM{b}_j^iF_{l,i}
  +\sum_{l=i+1}^{j-1}{\MM{b}_l^{i+1}}^{\CC{T}}\MM{b}_j^{i}F_{l,i}
  +\sum_{l=j+1}^{n}{\MM{b}_l^i}^{\CC{T}}\MM{b}_j^iF_{l,i}
  =L_i^j\MM{b}_i^j{\MM{e}_i}^{\CC{T}}-{\MM{e}_i}(L_i^j\MM{b}_i^j)^{\CC{T}},
\end{displaymath}   
where the matrix $L_i$ is the lower triangular part of 
$b_{1,2}F_{1,2}+\ldots+b_{i-1,n}F_{i-1,n}$ and we denote by
$L_{i}^{s}$, $s=i,j$, the matrix $L_i$ with zeros along its $s$-th row. 
Summing up with respect to $i$ and $j$, we obtain
\begin{Eqnarray*}
  \sum_{i=1}^{n-1}\sum_{j=i+1}^nb_{i,j}(-L_i^i\MM{b}_j^i){\MM{e}_i}^{\CC{T}}
  &=&\sum_{i=1}^{n-1}(-L_i^i)\left(\sum_{j=i+1}^nb_{i,j}\MM{b}_j^i\right) 
  {\MM{e}_i}^{\CC{T}}\\ 
  \sum_{i=1}^{n-1}\sum_{j=i+1}^nb_{i,j}L_i^j\MM{b}_i^j{\MM{e}_j}^{\CC{T}}
  &=&\sum_{j=2}^{n}\left(\sum_{i=1}^{j-1}b_{i,j}\MM{c}_i\right)
  {\MM{e}_j}^{\CC{T}} 
\end{Eqnarray*}
where we have used the notation $\MM{c}_i:=L_i\MM{b}_i^i$ for
$i=1,\dots,n-1.$  The cost of computing  the first sum is
$\frac43n^3$, while the cost of computing the second is
$2n^3$ operations.

Finally the terms
\begin{displaymath}
  \sum_{i=1}^{n-1}\sum_{j=i+1}^n\sum_{l=1}^ib_{i,j}b_{l,i}
  (\MM{b}_l^l\MM{e}_j^{\CC{T}}-\MM{e}_j{\MM{b}_l^l}^{\CC{T}})=
  \sum_{l=1}^{n-1}\MM{b}_l^l\left(\sum_{j=l+2}^n
    c_{j,l}\MM{e}_j^{\CC{T}}\right)-\left(\sum_{j=l+2}^n
    c_{j,l}\MM{e}_j\right){\MM{b}_l^l}^{\CC{T}}
\end{displaymath}
with $c_{j,l}=\sum_{i=l}^{j-1}b_{i,j}b_{l,i}$,
\begin{displaymath}
  \sum_{i=1}^{n-1}\sum_{j=i+1}^n\sum_{l=1}^ib_{i,j}b_{l,j}
  (\MM{b}_l\MM{e}_i^{\CC{T}}-\MM{e}_i\MM{b}_l^{\CC{T}})=
  \sum_{l=1}^{n-1}\MM{b}_l^l\left(\sum_{i=l+1}^{n-1}
    d_{i,l}\MM{e}_i^{\CC{T}}\right)-\left(\sum_{i=l+1}^{n-1}
    d_{i,l}\MM{e}_i\right){\MM{b}_l^l}^{\CC{T}}
\end{displaymath}
with $d_{i,l}=\sum_{j=i+1}^{n}b_{i,j}b_{l,j}$, and
\begin{Eqnarray*}
  &&\sum_{i=1}^{n-1}\sum_{j=i+1}^n\sum_{l=i+1}^{j-1}b_{i,j}b_{l,i}
  (\MM{b}_l^{i+1}\MM{e}_j^{\CC{T}}-\MM{e}_j{\MM{b}_l^{i+1}}^{\CC{T}})\\
  &=&\sum_{i=1}^{n-1}\sum_{l=i+1}^{n-1}b_{l,i}\left(\MM{b}_l^{i+1}
    {(\MM{b}_i^{l+1}-\MM{b}_i)}^{\CC{T}}-(\MM{b}_i^{l+1}-\MM{b}_i)
    {\MM{b}_l^{i+1}}^{\CC{T}}\right);
\end{Eqnarray*}
can be computed in about $\frac23n^3$, $\frac23n^3$ and $\frac12n^3$
operations respectively. Collecting the contributions of
all the terms in the sum we obtain a total count of  $7\frac12n^3$
operations. 

At the present time it is not clear that this method of computation of
$Q^2$ in the $\GG{so}(n)$ case is optimal form the point of view of
complexity theory. We did not try any other ordering of the basis 
elements and it is not at all certain that different orderings could give
better constants in front of the term $n^3$. 

Given that the construction of the (second-order) Strang splitting
carries a cost of $4n^3$ operations, the total flop count for
constructing a symmetric fourth order SKC approximation of an
exponential in $\GG{so}(n)$ by our algorithm is $11\frac12n^3$. This
is marginally better than obtaining an order-4
approximation by the Yo\u{s}ida technique from three Strang splittings
which, as pointed out in \cite{celledoni98atm}, requires $12n^3$
flops.

\setcounter{equation}{0}
\section{Sparse matrices}

In a naive formulation, the method of canonical coordinates of the
second kind is considerably too expensive for practical computation.
This, however, can be alleviated by the use of a sufficiently `sparse'
basis of the underlying Lie algebra $\GG{g}$. As explained in
Section~2, choosing a basis so that an overwhelming majority of
structure constants vanish renders the algorithm strikingly more
effective. It is important to emphasize that this has nothing to do
with the structure of the matrix $B\in\GG{g}$, which need not be
sparse. Yet, in most practical computations (in particular when $n$ is
large) one can expect $B$ to be sparse and structured. Good algorithms
should be able to exploit this phenomenon. 

In the case of SKC methods we identify two mechanisms that allow us to
exploit sparsity. Although this aspect of our methods is still a
matter for active investigation, the interim results are substantive
enough to warrant publication. For simplicity, we describe the first
mechanism just in the case of a tridiagonal $B\in\GG{so}(n)$, hence
\begin{displaymath}
  B=\left[
    \begin{array}{ccccc}
      0 & \beta_1 & 0 & \cdots & 0\\
      -\beta_1 & 0 & \ddots & \ddots & \vdots\\
      0 & \ddots & \ddots & \ddots & 0\\
      \vdots & \ddots & \ddots & 0 & \beta_{n-1}\\
      0 & \cdots & 0 & -\beta_{n-1} & 0
    \end{array}\right]=\sum_{k=1}^{n-1}\beta_k F_{k,k+1},
\end{displaymath}
where the matrices $F_{k,l}=\MM{e}_k\MM{e}_l^{\CC{T}}-\MM{e}_l\MM{e}_k
^{\CC{T}}$ have been introduced in Section~2. Since
\begin{displaymath}
  b_{k,l}=
  \begin{case}
    \beta_k, & l=k+1,\\
    0, & \mbox{otherwise,}
  \end{case}
\end{displaymath}
it is easy to substitute in the general formulae for the order-2
method:
\begin{displaymath}
  g_{i,j}'(0)=
  \begin{case}
    \beta_i, & \!\!\!\! j=i+1,\\
    0, & \!\!\!\! \mbox{otherwise,}
  \end{case}\quad  g_{i,j}''(0)=
  \begin{case}
    \beta_{i-1}\beta_i, & \!\!\!\! j=i-2,\\
    0, & \!\!\!\! \mbox{otherwise,}
  \end{case}
  \qquad 1\leq i<j\leq n,
\end{displaymath}
Arranging the elements of the basis in lexicographic order, we thus
obtain the second-order approximant
\begin{displaymath}
  \CC{e}^{\beta_{n-1}t F_{n-1,n}} \CC{e}^{\beta_{n-2}t F_{n-2,n-1}}
  \CC{e}^{\frac12\beta_{n-1}\beta_{n-1}t^2 F_{n-2,n}} \cdots
  \CC{e}^{\beta_2 tF_{2,3}} \CC{e}^{\frac12 \beta_2\beta_3 t^2F_{2,4}}
  \CC{e}^{\beta_1 tF_{1,2}} \CC{e}^{\frac12 \beta_1\beta_2 t^2F_{1,3}}.
\end{displaymath}
In other words, the cost of the approximation is just $\O{n}$ flops. 

Similar situation pertains to
\begin{displaymath}
  B=\left[
    \begin{array}{ccccc}
      \gamma_1 & \eta_1 & 0 & \cdots & 0\\
      \mu_1 & \gamma_2 & \ddots & & \vdots\\
      0 & \ddots & \ddots & \ddots & 0\\
      \vdots & & \ddots & \gamma_{n-1} & \eta_{n-1}\\
      0 & \cdots & 0 & \mu_{n-1} & \gamma_n
    \end{array}\right]\in \GG{sl}(n).
\end{displaymath}
Choosing the same basis and terminology as in Section 2 we can readily
ascertain that
\begin{Eqnarray*}
  g_{k,k-2}''(0)&=&\mu_{k-2}\mu_{k-1}\qquad k=3,4,\ldots,n\\
  g_{k,k-1}''(0)&=&-(\gamma_{k-2}-2\gamma_{k-1}+\gamma_k)\mu_{k-1},
  \qquad k=3,4,\ldots,n,\\
  g_{k,k+1}''(0)&=&(\gamma_{k-1}-2\gamma_k+\gamma_{k+1})\eta_k,\qquad
  k=2,3,\ldots, n-1,\\
  g_{k,k+2}''(0)&=&-\eta_k\eta_{k+1},\qquad k=1,2,\ldots,n-2,\\
  g_{k,l}''(0)&=&0,\qquad |k-l|\geq 3
\end{Eqnarray*}
and $g_k''(0)=-\eta_k\mu_k$, $k=1,2,\ldots,n-1$. Thus, a second-order
approximant to a tridiagonal $B\in\GG{sl}(n)$ is itself quindiagonal
and its computation requires just $\O{n}$ flops.

Higher-order approximants and matrices with greater bandwidth lend
themselves to similar treatment, although the savings are less
striking. In a sense, the situation is parallel to that of
approximating $\exp(tB)$ by a rational approximant, when savings
accrue from sparse matrix-inversion methods, except that in our case
the result is assured to belong to the right Lie group.

Another observation which is highly pertinent to the approximation of
exponentials of sparse matrices has been made in
\cite{iserles99hli}. Suppose that $B$ is a banded matrix of bandwidth
$s\geq 3$. In general, $F(t)=\exp(tB)$ is a dense matrix. Yet, as is
easy to illustrate by computer experiments, $F(t)$ is very near to a
banded matrix. Specifically, given $\varepsilon>0$, there exists
$r=r(t,\varepsilon)\geq s$ such that all the elements of $F(t)$
outside a band of width $r$ are less than $\varepsilon$ in magnitude.
Moreover, tight upper bounds on $r$ can be derived with relative
ease. The idea thus is to set to zero all the elements outside
bandwidth $r$. The outcome is a banded approximant to the
exponential. Moreover, with an appropriate choice of basis elements,
this means that the functions $\alpha_i$ are {\em set to zero\/} for
elements that possess terms exclusively outside the band. Consequently
corresponding exponentials equal identity and need not be included in
the product. Thus, the cost scales with the size $r$ of the
bandwidth. Similar phenomenon has been already encountered in the context
of $\GG{so}(n)$ and $\GG{sl}(n)$, when our choice of basis and order
has implied a banded structure of the exponential. The present
mechanism is different, even if the net outcome is similar. 

\setcounter{equation}{0}
\section{Numerical experiments}

Our numerical experiments are organized as follows. We fist consider
a test on random matrices  in $\GG{so}(50)$, illustrating the
performance of methods based on the use of second kind coordinates
techniques for full and sparse matrices. The third and last example is
the solution of a third-order ODE using {\rm
  Runge--Kutta/Munthe-Kaas\/} (RK/MK) methods described in
\cite{munthe-kaas97hor}.  We use the {\tt Matlab} toolbox {\tt
  DiffMan\/} for the integration of ODEs on manifolds, comparing
the  usual implementation of RK/MK methods, whereby the the
exponential is approximated to machine accuracy, with a version of the
methods obtained using the time-symmetric fourth-order approximation
from Section~4.

All experiments have been performed in {\tt Matlab} and we have
computed the error while comparing the results with the built-in
function {\tt expm} which calculates the exponential to nearly machine
accuracy.   

We evaluated the 
the error computing $\|\CC{e}^{-tB}F(tB)-I\|_{\rm F}$ where $F(tB)$ is the
SKC approximation of $\exp(tB)$ and $\|\cdot\|_{\rm F}$ denotes the
Frobenius norm.  
The matrices have been generated randomly using the {\tt Matlab} function
{\tt rand} and scaling the Frobenius norm so that
$\|B\|_{\rm F}=1$. 

We approximate $\exp(tB)$ with a single step of the methods for different
values of $t$, ($t=1/2^k$ and $k=1,\ldots ,5$).

\begin{figure}[Htbp]
  \begin{center}
    \leavevmode
    \psfig{file=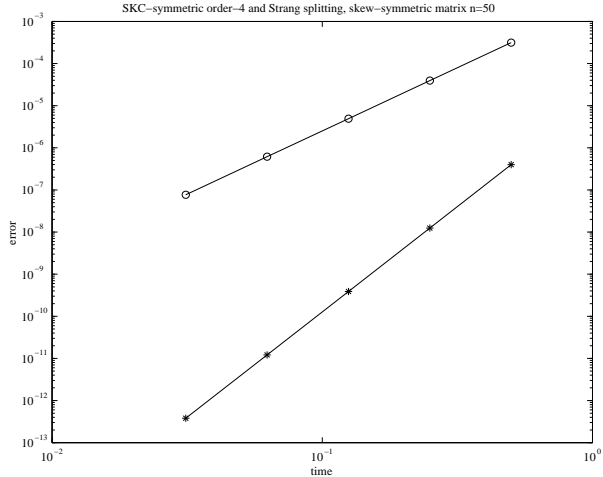,width=8cm}
    \caption{Error versus time in the $\GG{so}(50)$ (full case).}
    \label{fig:2}
  \end{center}
\end{figure}

In both the first two figures the norm of the error is plotted (along
the $y$-axis) to a logarithmic scale with respect to $t$.
Figure~\ref{fig:2} reports the results of our first test, where we
have considered a full matrix in $\GG{so}(50)$.  In the plots the
error norm is indicated with the symbols `$*$' (SKC, time symmetric,
order 4) and `$\circ$' (Strang splitting, order 2).

\begin{figure}[Htbp]
  \begin{center}
    \leavevmode
    \psfig{file=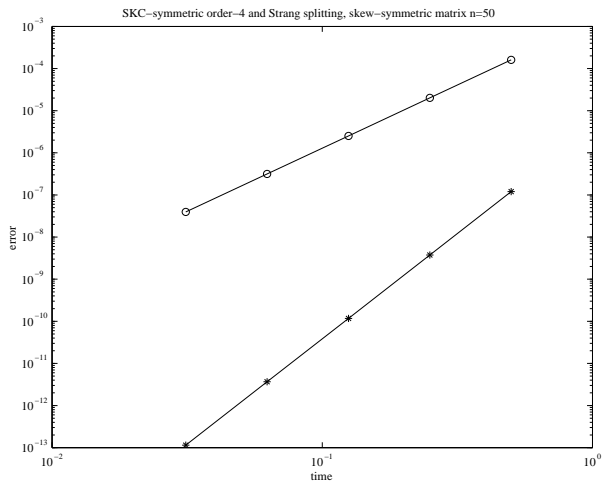,width=8cm}
    \caption{Error versus time in the $\GG{so}(50)$ sparse case.}
    \label{fig:3}
  \end{center}
\end{figure}

In the next example, illustrated in Figure~\ref{fig:3}, the same
methods have been applied to a sparse matrix in $\GG{so}(50)$, with
four non-zero diagonals (i.e., bandwidth 5). In both the examples the
methods give the correct order.  In the second case, however, the
count of flops is drastically reduced. We counted the number of flops
using the {\tt Matlab} function {\tt flops}. In the first case the
cost for constructing $Q^2$ amounts to $9.62n^3$ while in the second
we counted $0.95n^3$ flops. As it is easy to understand, the described
implementation of the methods allows to take advantage immediately of
the sparsity structure of the matrix $B$, working directly on the
nonzeros entries of $B$ {\em \`a la\/} Section~4.

The last example is concerned with the use of the techniques described in
this paper in substituting the exponentials computed to machine
accuracy in the integration methods of \cite{munthe-kaas97hor}. 
The experiments have been performed using the {\tt Matlab} toolbox
{\tt DiffMan.\/} We use a RK/MK method of order four. 

The example is a problem whose solution is the soliton originating in
the {\em Korteweg--de Vries\/} (KdV) equation. It is a third-order ODE
obtained performing a symmetry reduction on the KdV equation.
The resulting ODE can be written as a three-dimensional system,
\begin{displaymath}
  y'=\left[
    \begin{array}{ccc}
      0&1 &0\\
      0&0&1\\
      -9y(2)&3&0\\
    \end{array}
  \right]y
\end{displaymath}
with $y(0)=[1,0,-1.5]^{\CC{T}}$ and $t\in [0,5]$. The solution of the
ODE $f=y_1(t)$ can be easily derived explicitely and it is
$f(t)=\alpha\CC{sech}\,(t\beta)$, $\alpha=1$, $\beta=1/2\sqrt(3).$

\begin{figure}[Htbp]
  \begin{center}
    \leavevmode
    \psfig{file=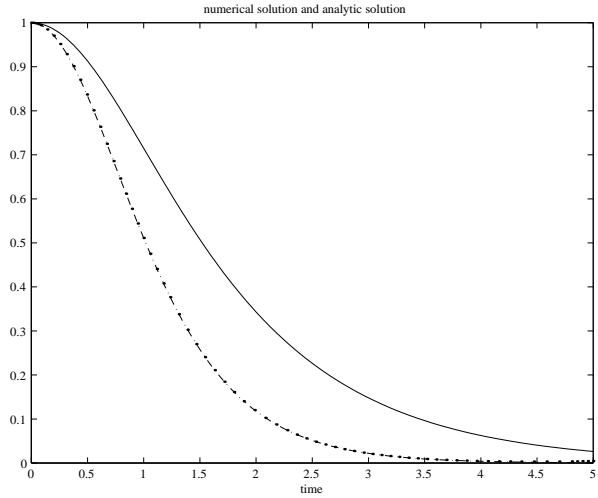,width=8cm}
    \caption{The soliton originating in the KdV equation.}
    \label{fig:4}
  \end{center}
\end{figure}

In Figure \ref{fig:4} we plot the analytic solution (solid line) on a
grid of $161$ points. The dotted line is the numerical solution
obtained with the {\tt Matlab} routine {\tt ode45} with absolute and
relative tolerance $1.0e-4$. The method produced this solution in $69$
steps, and it was implemented with step-size control procedure.  The
dashed-dotted line is the numerical solution obtained with the RK/MK
method using SKC symmetric tecnhniques for the approximation of the
exponential, with fixed step-size $h$.

\begin{figure}[Htbp]
  \begin{center}
    \leavevmode
    \psfig{file=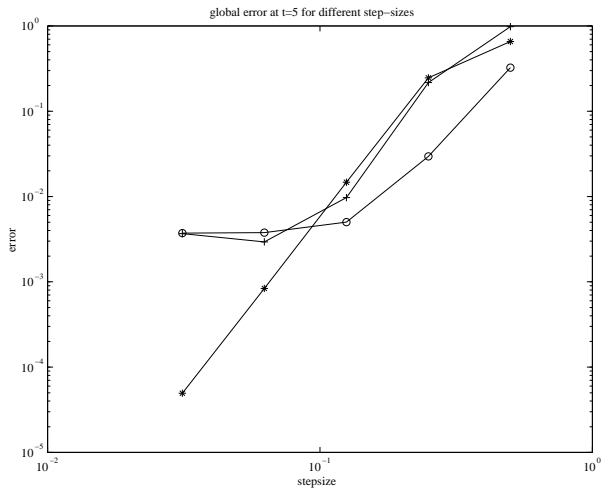,width=8cm}
    \caption{RK/MK: global error at $t=5$ with {\tt expm} and SKC}
    \label{fig:5}
  \end{center}
\end{figure}

In Figure \ref{fig:5} we plot the error (along the $y$-axis) with
respect to the numerical solution obtained with the {\tt Matlab}
routine {\tt ode45} to a logarithmic scale, versus the stepsize
$h=1/2^k$ and $k=1,\ldots ,5$ for the cases of the implementation of
RK/MK with the {\tt expm} function of {\tt Matlab} (marked with $+$)
and approximating $\exp$ to order four with a SKC technique ($\circ$).
The line marked with $\ast$ representes the error of the numerical
solution given by the RK/MK method implemented with SKC technique for
the approximation of the exponential, measured with respect to the
numerical solution obtained by the same method with the use of the
exact exponential ({\tt expm} routine of {\tt Matlab}).

It is interesting to note in this case that substituting the exact
exponential with suitable fourth-order approximant does not lead
to a significant deterioration in the quality of the RK/MK method and
the overall error does not change much.
Note that in the present case the primary variable is a vector, rather
than a matrix. In general, if the underlying ODE  can be written in a
vector form, i.e.\ as an action of a Lie group on $\BB{R}^n$, we need
to approximate $\exp(tB)\MM{v}$, where $\MM{v}\in\BB{R}^n$, rather
than the matrix $\exp(tB)$. This leads to obvious savings in the SKC
techniques, similarly, say, to the approach of rational functions. In
particular, the cost of composing exponentials is $\O{n^2}$, rather
than $\O{n^3}$, operations.

\section*{Acknowledgments}

The authors are grateful to Brynjulf Owren for many fruitful
discussions, to Per Christian Moan for bringing the reference
\cite{wei63ogr} to their attention, to the Numerical Analysis group of
DAMTP Cambridge, and to the Geometric Integration members during the
fall semester 1998 at MSRI Berkeley. Research at MSRI is supported in
part by NSF grant DMS-9701755. 

\bibliographystyle{agsm}
\bibliography{geom_int,hamilt,elena}

\appendix
\setcounter{equation}{0}
\section{Appendix}

For completeness, we present a proof of Lemma 1. Note that a
comprehensive treatment of this subject matter, inclusive of the
non-symmetric case, has been presented in a different context by
\citeasnoun{chacon91rff}. For our purposes, however, it is sufficient
to derive the first  term of the expansion.

\setcounter{theorem}{0}
\setcounter{section}{4}
\setcounter{equation}{5}
\begin{lemma}
  The leading error term in the Strang splitting is
  \begin{equation}
    \label{eq:A.1}
    \Frac{1}{12}\sum_{l=2}^s [C_1+\cdots+C_{l-1}+\Frac12
    C_l,[C_1+\cdots+C_{l-1},C_l]]. 
  \end{equation}
\end{lemma}

\begin{proof}
  Letting
  \begin{Eqnarray*}
    F_1(t)&=&\CC{e}^{tC_1},\\
    F_l(t)&=&\CC{e}^{tC_l/2}F_{l-1}\CC{e}^{tC_l/2},\qquad l=2,3,\ldots,s,
  \end{Eqnarray*}
  we can verify at once that $F_s$ is precisely the Strang
  splitting. We assume that
  \begin{displaymath}
    F_l(t)=\exp[t(C_1+\cdots+C_l)+\Frac{1}{12} Q_l t^3+\O{t^4}],\qquad
    l=1,2,\ldots,s.
  \end{displaymath}
  We use the BCH formula:
  \begin{Eqnarray*}
    F_l(t)&=&\exp (\Frac12 tC_l) \exp[t(C_1+\cdots+C_{l-1})
    +\Frac{1}{12} Q_{l-1}t^3+\O{t^4}] \exp (\Frac12 tC_l)\\
    &=&\exp \{t(C_1+\cdots+C_{l-1}+\Frac12 C_l) +\Frac14 t^2 [C_l,
    C_1+\cdots +C_{l-1}]\\
    &&\mbox{}\quad+\Frac{1}{24} t^3 [\Frac12 C_l-(C_1+\cdots
    +C_{l-1}), [C_l,C_1+\cdots+C_{l-1}]]\\
    &&\mbox{}\quad+\Frac{1}{12}t^3Q_{l-1}+\O{t^4}\} \exp (\Frac12
    tC_l)\\ 
    &=&\exp \{ t(C_1+\cdots+C_l)+\Frac14 t^2 [C_l,C_1+\cdots+C_{l-1}]\\
    &&\mbox{}\quad+\Frac{1}{24} t^3 [\Frac12 C_l-(C_1+\cdots +C_{l-1}),
    [C_l,C_1+\cdots+C_{l-1}]]\\
    &&\mbox{}\quad+\Frac14 t^2 [C_1+\cdots
    +C_{l-1}+\Frac12 C_l,C_l] +\Frac{1}{16} t^3
    [[C_l,C_1+\cdots+C_{l-1}], C_l]\\
    &&\mbox{}\quad +\Frac{1}{24} t^3
    [C_1+\cdots+C_{l-1}, [C_1+\cdots+C_{l-1}+\Frac12 C_l,C_l]]
    +\Frac{1}{12}t^3Q_{l-1}+\O{t^4}\}. 
  \end{Eqnarray*}
  However,
  \begin{displaymath}
    [C_l,C_1+\cdots+C_{l-1}]+[C_1+\cdots+C_{l-1}+\Frac12 C_l,C_l]=O,
  \end{displaymath}
  thereby annihilating the $t^2$ term, and
  \begin{Eqnarray*}
    &&\Frac{1}{24}[\Frac12 C_l-(C_1+\cdots +C_{l-1}), [C_l,C_1+\cdots
    +C_{l-1}]]+\Frac{1}{16}[[C_l,C_1+\cdots+C_{l-1}], C_l]\\
    &&\mbox{} +\Frac{1}{24}[C_1+\cdots+C_{l-1}, [C_1+\cdots+C_{l-1}
    +\Frac12 C_l,C_l]] \\
    &=&\Frac{1}{12} [C_1+\cdots+C_{l-1}+\Frac12
    C_l,[C_1+\cdots+C_{l-1}, C_l]].
  \end{Eqnarray*}
  Therefore
  \begin{displaymath}
    Q_l=[C_1+\cdots+C_{l-1}+\Frac12 C_l,[C_1+\cdots+C_{l-1}, C_l]]
    +Q_{l-1}. 
  \end{displaymath}
  Since $Q_1=O$, the expression \R{eq:A.1} follows by summing the above
  formula for $l=2,3,\ldots,s$.
\end{proof}

\end{document}